\documentclass[reqno,12pt]{amsart}
\usepackage{a4wide}

\usepackage{amssymb,graphicx}

\newcommand{\CH}[1][2]{\ensuremath{\mathbb{C}H^{#1}}}
\newcommand{\C}{\ensuremath{\mathbb{C}}}
\newcommand{\R}{\ensuremath{\mathbb{R}}}

\newcommand{\sech}{\mathop{\mbox{sech}}}

\newcommand{\enabla}{\ensuremath{\bar{\nabla}}}
\newcommand{\eR}{\ensuremath{\bar{R}}}
\newcommand{\tr}{\mathop{\mbox{tr}}}

\newtheorem{theorem}{Theorem}[section]
\newtheorem{proposition}[theorem]{Proposition}
\newtheorem{lemma}[theorem]{Lemma}
\newtheorem{corollary}[theorem]{Corollary}

\usepackage{color}

\begin{document}

\title[Real hypersurfaces with constant principal curvatures]
{\Large Real hypersurfaces\\ with constant principal curvatures\\
 in the complex hyperbolic plane}

\author{J{\"u}rgen Berndt, Jos{\'e} Carlos D{\'\i}az-Ramos}

\maketitle

\section{Introduction}

\'{E}lie Cartan classified in \cite{Ca38} all connected hypersurfaces $M$ with
constant principal curvatures in the real hyperbolic space $\R H^n$, $n \geq 3$.
His classification exhibits two remarkable features. Firstly, the number $g$
of distinct principal curvatures has an upper bound independent of the dimension
$n$. In fact, Cartan showed that $g \leq 2$. Secondly,
every connected real hypersurface with constant principal curvatures
in $\R H^n$ is an open part of a homogeneous hypersurface. Therefore,
assuming $M$ is complete, the constancy of principal curvatures
is equivalent to the existence of a closed subgroup of the isometry group $G$ of $\R H^n$
such that $M$ is an orbit of $G$.

We are interested in the corresponding questions in the complex hyperbolic space
$\C H^n$, $n \geq 2$. We summarize briefly the known facts so far.
Cartan's argument for $g \leq 2$ in $\R H^n$ relies on the Gau\ss-Codazzi equations.
The structure of these equations in $\C H^n$ is too complicated in general, but
simplifies considerably for Hopf hypersurfaces. A real hypersurface $M$ in $\C H^n$
is a Hopf hypersurface if the Hopf foliation on $M$ is totally geodesic. The Hopf
foliation on $M$ consists of the leaves of the one-dimensional distribution on $M$
that is obtained by rotating the normal bundle of $M$ into the tangent bundle of $M$
by means of the complex structure $J$ of $\C H^n$. If $\xi$ is a unit normal vector field
on $M$, then $M$ is a Hopf hypersurface if and only if $J\xi_p$ is a principal curvature
vector of $M$ at each point $p \in M$. The first author classified in \cite{Be89} all
connected Hopf hypersurfaces with constant principal curvatures
in $\C H^n$, $n \geq 2$. It turns out that
$g \in \{2,3\}$ and that $M$ is an open part of a homogeneous real hypersurface in $\C H^n$.

For some time it was believed that every homogeneous real hypersurface in $\C H^n$ is
a Hopf hypersurface, a fact known to be true in the complex projective space $\C P^n$
according to Takagi's classification in \cite{Ta73}. Surprisingly, Lohnherr \cite{Lo98}
constructed in his PhD thesis an example of a homogeneous real hypersurface $W^{2n-1}$
in $\C H^n$ which is not a Hopf hypersurface. Consider a horocycle $H$ in a totally
geodesic $\R H^2 \subset \C H^2 \subset \C H^n$, and attach to each point $p \in H$
the totally geodesic $\C H^{n-1} \subset \C H^n$ which is perpendicular to $H$ at $p$.
The resulting ruled real hypersurface $W^{2n-1}$ is a minimal homogeneous real hypersurface
in $\C H^n$. An alternative Lie theoretic construction of $W^{2n-1}$ has been presented by the
first author in \cite{Be98}.

The first author and Br\"{u}ck constructed in \cite{BB01} more examples of
homogeneous real hypersurfaces in $\C H^n$ which are not Hopf hypersurfaces.
Recently, the first author and Tamaru classified in \cite{BT06} the cohomogeneity one actions
on $\C H^n$ up to orbit equivalence. This of course provides the classification of homogeneous real
hypersurfaces in $\C H^n$, and confirms that there are no further homogeneous
real hypersurfaces in $\C H^n$ apart from the known ones.
Any such hypersurface has constant principal curvatures with
$g \in \{2,3,4,5\}$. The authors classified in \cite{BD06} all connected real hypersurfaces
in $\C H^n$, $n \geq 3$, with at most three distinct constant principal curvatures.
A consequence of this classification is that any such hypersurface is an open part of a homogeneous
real hypersurface in $\C H^n$. The methods developed in \cite{BD06} do not work for
the case $n = 2$. The purpose of this paper is to settle this remaining case by a different method.

\begin{theorem}\label{mainresult}
Let $M$ be a connected real hypersurface in the complex hyperbolic
plane $\C H^2$. Then $M$ has constant principal curvatures if and
only if it is is an open part of a homogeneous hypersurface in $\C
H^2$, that is, $M$ is an open part of
\begin{itemize}
\item[(i)] a geodesic hypersphere in $\C H^2$, or

\item[(ii)] a horosphere in $\C H^2$, or

\item[(iii)] a tube around a totally geodesic $\CH[1]\subset\CH$,
or

\item[(iv)] a tube around a totally geodesic $\R H^2\subset\CH$,
or

\item[(v)] a ruled minimal real hypersurface $W^{3}\subset\CH$, or
one of its equidistant hypersurfaces.
\end{itemize}
\end{theorem}

The second author has been supported by grants from
Ministerio de Ciencia y Tecnolog\'{i}a (BFM 2003-02949)
and Xunta de Galicia (PGIDIT 04PXIC20701PN).

\section{Preliminaries}\label{scFormulae}

Let $\CH$ be the complex hyperbolic plane equipped with the Fubini
Study metric $\langle \cdot , \cdot \rangle$ of constant
holomorphic sectional curvature $-1$. We denote by $\bar{\nabla}$
and $\bar{R}$ the Levi Civita covariant derivative and the
Riemannian curvature tensor of $\CH$, respectively, using the sign
convention $\bar{R}_{XY} = [\bar{\nabla}_X,\bar{\nabla}_Y] -
\bar{\nabla}_{[X,Y]}$. Then
\[
\bar{R}_{XY}Z = -\frac{1}{4}\Big(\langle Y,Z\rangle X-\langle
X,Z\rangle Y +\langle JY,Z\rangle JX-\langle JX,Z\rangle JY
-2\langle JX,Y\rangle JZ\Big),
\]
where $J$ is the complex structure of $\CH$. We also write
$\bar{R}_{XYZW}=\langle \bar{R}_{XY}Z,W\rangle$.

Let $M$ be a connected real hypersurface of $\CH$. We denote by
$\nabla$ and $R$ the Levi Civita covariant derivative and the
Riemannian curvature tensor of $M$, respectively. By $TM$ and $\nu
M$ we denote the tangent bundle and the normal bundle of $M$, and
by $\Gamma(TM)$ and $\Gamma(\nu M)$ we denote the module of all
smooth vector fields tangent and normal to $M$, respectively. Let
$X,Y,Z,W \in \Gamma(TM)$ and $\xi\in \Gamma(\nu M)$ be a
(local) unit normal vector field on $M$.

The Levi Civita covariant derivatives of $M$ and $\CH$ are related
by the Gau\ss\ formula
\[
\enabla_X Y=\nabla_X Y+\langle S X,Y \rangle \xi,
\]
where $S$ is the shape operator of $M$ with respect to $\xi$. The
Weingarten formula is
\[
\enabla_X\xi = -S X.
\]
The fundamental equations of second order of interest to us are
the Gau\ss\ equation
\[
\eR_{XYZW}=R_{XYZW} - \langle SY,Z\rangle\langle SX,W \rangle
    +\langle SX,Z\rangle\langle SY,W\rangle
\]
and the Codazzi equation
\[
\eR_{XYZ\xi}=\langle (\nabla_XS)Y-(\nabla_YS)X,Z \rangle.
\]

We assume from now on that $M$ has constant principal curvatures.
For each principal curvature $\lambda$ of $M$ we denote by
$T_\lambda$ the distribution on $M$ formed by the principal
curvature spaces of $\lambda$. By $\Gamma(T_\lambda)$ we denote
the set of all smooth sections in $T_\lambda$, that is, all smooth
vector fields on $M$ satisfying $SX = \lambda X$. The Codazzi
equation readily implies
\begin{lemma}\label{thThreeEigenv} For all $X\in\Gamma(T_{\lambda_i})$,
$Y\in\Gamma(T_{\lambda_j})$ and $Z\in\Gamma(T_{\lambda_k})$ we
have
$$
\bar{R}_{XYZ\xi} = ({\lambda_j}-{\lambda_k})\langle\nabla_X
Y,Z\rangle
    -(\lambda_i-{\lambda_k})\langle\nabla_Y X,Z\rangle.
$$
\end{lemma}

Putting $\lambda_i = \lambda_k$ in Lemma \ref{thThreeEigenv} and
then interchanging $Y$ and $Z$ yields

\begin{lemma}\label{thTwoEigenv}
For all $X,Y \in \Gamma(T_{\lambda_i})$ and $Z \in
\Gamma(T_{\lambda_j})$ with $\lambda_i \neq {\lambda_j}$ we have
\[
4({\lambda_j}-{\lambda_i})\langle\nabla_XY,Z\rangle = \langle JY,Z
\rangle \langle X,J\xi\rangle + \langle JX,Y\rangle \langle
Z,J\xi\rangle + 2\langle JX,Z\rangle \langle Y,J\xi \rangle.
\]
\end{lemma}

\section{The Classification}\label{scPrincipalCurvatures}

Let $M$ be a connected real hypersurface of $\CH$ with constant
principal curvatures. Since our classification problem is of local
nature, we can assume that $M$ is orientable. Let $\xi$ be a
global unit normal vector field on $M$ and denote by $g$ the
number of distinct constant principal curvatures of $M$. As an
immediate consequence of Lemma \ref{thThreeEigenv} we have

\begin{lemma}\label{thNoUmbilical}
$g \geq 2$.
\end{lemma}

\begin{proof}
Suppose that $g=1$. Then Lemma \ref{thThreeEigenv} implies
$\bar{R}_{XYZ\xi}=0$ for all $X,Y,Z \in \Gamma(TM)$. In
particular, $0=4\eR_{J\xi YZ\xi} =\langle JY,Z\rangle$ for all
$Y,Z\in\Gamma(TM)$, which means that $M$ is a totally real
submanifold of $\C H^2$, and hence $\dim M \leq 2$, which is a
contradiction to $\dim M = 3$.
\end{proof}

\begin{proposition}\label{th2Constant}
If $g = 2$ then $M$ is an open part of
\begin{itemize}
\item[(i)] a geodesic hypersphere in $\C H^2$, or

\item[(ii)] a horosphere in $\C H^2$, or

\item[(iii)] a tube around a totally geodesic $\CH[1]\subset\CH$,
or

\item[(iv)] a tube of radius $r = \ln(2+\sqrt3)$ around a totally
geodesic $\R H^2\subset\CH$.

\end{itemize}
\end{proposition}

\begin{proof}
We denote by $\lambda_i$ the principal curvature of $M$ with
multiplicity $i \in \{1,2\}$. Suppose there is a point $p \in M$
such that $J\xi_p$ is not a principal curvature vector. Then there
exists an open neighborhood $N$ of $p$ in $M$ such that
$J\xi=b_1U_1+b_2U_2$ for some unit vector fields
$U_i\in\Gamma(T_{\lambda_i})$ on $N$ and smooth functions
$b_i:N\to\R$ with $b_i(q) \neq 0$ for all $q \in N$. As the rank
of the distribution $T_{\lambda_2}$ is $2$, there exists a nonzero
vector field $V\in\Gamma(T_{\lambda_2})$ on $N$ such that $\langle
V,U_2\rangle=0$. By construction we have $\langle
JV,\xi\rangle=-\langle V,J\xi \rangle = 0$ and $\langle
JV,V\rangle=0$. From Lemma \ref{thThreeEigenv} we get
$0=4\eR_{U_2VU_2\xi}=3b_2\langle JV,U_2\rangle$. Thus we have
$\langle JV,U_2\rangle=0$ on $N$, and hence $b_1 \langle
JV,U_1\rangle= - \langle V, b_1JU_1 \rangle = \langle V,\xi +
b_2JU_2\rangle =0$. This implies $\langle JV,U_1\rangle =0$
on $N$ as well. Altogether this shows that $V = 0$ on $N$, which
is a contradiction. We conclude that $J\xi_p$ is a principal
curvature vector of $M$ at each point $p \in M$. The proposition
then follows from \cite{Be89}.
\end{proof}

In view of the above results we can assume $g=3$ from now on. The
real hypersurfaces with constant principal curvatures in $\C H^n$,
$n \geq 3$, have been classified by the authors in \cite{BD06}.
Unfortunately, the proof in \cite{BD06} does not work for $n=2$.
For this reason we need to develop here a different method for $n
= 2$.

We denote by $\lambda_1$, $\lambda_2$ and $\lambda_3$ the three
distinct principal curvatures of $M$. As $\dim M = 3$, each of
these principal curvatures has multiplicity one. Let $U_1,U_2,U_3$
be a local orthonormal frame field on $M$ with $U_i \in
\Gamma(T_{\lambda_i})$, that is $SU_i = \lambda_iU_i$. Since we
are interested only in the local structure of $M$, we can assume
without loss of generality that $U_1,U_2,U_3$ is a global
orthonormal frame field on $M$. Then we can write
$$J\xi = b_1U_1 + b_2U_2 + b_3U_3$$
with the smooth functions $b_i : M \to \R$ defined by $b_i =
\langle J\xi , U_i \rangle$, $i \in \{1,2,3\}$. Note that $$ b_1^2
+ b_2^2 + b_3^2 = 1.$$ In the following we assume that all indices
are taken modulo 3. As
\[ JU_i  =
\langle JU_i,U_{i+1} \rangle U_{i+1} - \langle JU_{i+2}, U_i \rangle
U_{i+2} - b_i\xi,
\]
we get
\[
0=\langle U_i,\xi\rangle=\langle JU_i,J\xi\rangle=
\langle JU_i,U_{i+1} \rangle b_{i+1} - \langle JU_{i+2}, U_i \rangle b_{i+2}
\]
This gives a system of three linear equations, and we easily see that the vector
$(b_1,b_2,b_3)$ is in the real span of $(\langle
JU_2,U_3\rangle,\langle JU_3,U_1\rangle,\langle JU_1,U_2\rangle)$.
From $b_1^2+b_2^2+b_3^2=1$ we get
\[
3 = \sum_{i=1}^3\langle U_i,U_i\rangle=\sum_{i=1}^3\langle
JU_i,JU_i\rangle=2\Bigl(\langle JU_2,U_3\rangle^2+\langle
JU_3,U_1\rangle^2+\langle JU_1,U_2\rangle^2\Bigr)+1.
\]
Thus $(\langle JU_2,U_3\rangle,\langle JU_3,U_1\rangle, \langle
JU_1,U_2\rangle)$ is a unit vector in $\R^3$, and hence we must have
$(b_1,b_2,b_3) = \pm (\langle JU_2,U_3\rangle,\langle JU_3,U_1\rangle, \langle
JU_1,U_2\rangle)$. Without loss of generality we can assume that
we have the plus sign (otherwise replace $\xi$ by $-\xi$). Thus we have proved
\begin{equation}\label{thJUixi}
\langle JU_i,U_{i+1}\rangle =b_{i+2}.
\end{equation}

We now introduce the following notation:
$$ d_i = \lambda_{i+1} - \lambda_{i+2}\ ,\ x_i = \langle
\nabla_{U_i}U_{i+1},U_{i+2} \rangle.$$
Putting $X = U_{i+1}$, $Y = U_{i+2}$ and $Z = U_i$ in Lemma \ref{thThreeEigenv},
and using (\ref{thJUixi}) and $b_1^2+b_2^2+b_3^2 = 1$, we obtain
\begin{equation}\label{thbi}
3b_i^2 = 1 + 4d_{i+2}x_{i+2} - 4d_{i+1}x_{i+1}.
\end{equation}
As $\nabla_{U_i}U_i=
\langle\nabla_{U_i}U_i,U_{i+1}\rangle U_{i+1}
+\langle\nabla_{U_i}U_i,U_{i+2}\rangle U_{i+2}$,
Lemma \ref{thTwoEigenv} and (\ref{thJUixi}) imply
\begin{equation} \label{thNablas}
4\nabla_{U_i}U_i =
    -3b_ib_{i+2}d_{i+2}^{-1}U_{i+1}-3b_ib_{i+1}d_{i+1}^{-1}U_{i+2}.
\end{equation}

We now calculate the differential $db_i$ of the function $b_i = \langle U_i,J\xi\rangle$.
We have $\bar{\nabla}_{U_i}J\xi = J\bar{\nabla}_{U_i}\xi = -\lambda_iJU_i$
and hence $\langle U_{i},\enabla_{U_i}J\xi\rangle =0$. This implies
$db_i(U_i) = b_{i+1}\langle\nabla_{U_i}U_{i},U_{i+1}\rangle
+ b_{i+2}\langle\nabla_{U_i}U_{i},U_{i+2}\rangle$, and using (\ref{thNablas})
we get
\begin{equation}\label{thdbii}
4db_i(U_i) = 3b_1b_2b_3d_id_{i+1}^{-1}d_{i+2}^{-1}.
\end{equation}
Using (\ref{thJUixi}) and (\ref{thNablas}) we get from
\begin{eqnarray*} db_i(U_{i+1}) &=& b_{i+1}\langle\nabla_{U_{i+1}}U_i,U_{i+1}\rangle
+ b_{i+2}\langle\nabla_{U_{i+1}}U_i,U_{i+2}\rangle - \lambda_{i+1}
\langle U_i,JU_{i+1} \rangle \\ &=&
-b_{i+1}\langle\nabla_{U_{i+1}}U_{i+1},U_i\rangle
- b_{i+2}\langle\nabla_{U_{i+1}}U_{i+2},U_i\rangle + \lambda_{i+1}
\langle JU_i,U_{i+1} \rangle
\end{eqnarray*}
that
\begin{equation}\label{thdbii1}
4db_i(U_{i+1}) = b_{i+2}(3b_{i+1}^2d_{i+2}^{-1} + 4\lambda_{i+1} - 4x_{i+1}).
\end{equation}
In a similar way we obtain
\begin{equation}\label{thdbii2}
4db_i(U_{i+2}) = b_{i+1}(3b_{i+2}^2d_{i+1}^{-1} - 4\lambda_{i+2} + 4x_{i+2}).
\end{equation}

{From} the Gau\ss\ equation and (\ref{thJUixi}) we get
\[
4R_{U_iU_{i+1}U_iU_{i+1}} = 1-4\lambda_i\lambda_{i+1}+3b_{i+2}^2.
\]
On the other hand, by definition,
\[
R_{U_iU_{i+1}U_iU_{i+1}}=
\langle\nabla_{U_i}\nabla_{U_{i+1}}U_i,U_{i+1}\rangle
-\langle\nabla_{U_{i+1}}\nabla_{U_{i}}U_i,U_{i+1}\rangle
-\langle\nabla_{[U_i,U_{i+1}]}U_i,U_{i+1}\rangle.
\]
The first term on the right-hand side can be calculated using
(\ref{thNablas}), (\ref{thdbii1}) and (\ref{thdbii2}) as follows
\[
\everymath{\displaystyle}
\begin{array}{rcl}
&&16\langle\nabla_{U_i}\nabla_{U_{i+1}}U_i,U_{i+1}\rangle\\[2ex]
&=&\langle\nabla_{U_i}(12b_{i+1}b_{i+2}d_{i+2}^{-1}U_{i+1}-16x_{i+1}U_{i+2}),
    U_{i+1}\rangle\\[2ex]
&=&12(db_{i+1}(U_i)b_{i+2}+b_{i+1}db_{i+2}(U_i))d_{i+2}^{-1}+16x_ix_{i+1}.\\[2ex]
&=&9(b_{i+1}^2d_{i+1}^{-1} + b_{i+2}^2d_{i+2}^{-1})b_i^2d_{i+2}^{-1}
+ 12(x_i-\lambda_i)(b_{i+2}^2-b_{i+1}^2)d_{i+2}^{-1} + 16x_ix_{i+1}.
\end{array}
\]
The second and third term can be calculated in a similar way:
\[
\everymath{\displaystyle}
\begin{array}{lcl}
&&16\langle\nabla_{U_{i+1}}\nabla_{U_{i}}U_i,U_{i+1}\rangle\\[2ex]
&=&-12\langle\nabla_{U_{i+1}}(b_ib_{i+2}d_{i+2}^{-1}U_{i+1}-
     b_ib_{i+1}d_{i+1}^{-1}U_{i+2}, U_{i+1}\rangle\\[2ex]
&=&-12(db_i(U_{i+1})b_{i+2}+b_idb_{i+2}(U_{i+1}))d_{i+2}^{-1}-
    9b_i^2b_{i+1}^2d_i^{-1}d_{i+1}^{-1}\\[2ex]
    &=&-9(b_i^2d_i^{-1} + b_{i+2}^2d_{i+2}^{-1})b_{i+1}^2d_{i+2}^{-1}
+ 12(x_{i+1}-\lambda_{i+1})(b_{i+2}^2-b_i^2)d_{i+2}^{-1} - 9b_i^2b_{i+1}^2d_i^{-1}d_{i+1}^{-1},\\[2ex]
&&16\langle\nabla_{[U_i,U_{i+1}]}U_i,U_{i+1}\rangle \\[2ex]
&=&-16\langle\nabla_{U_i}U_i,U_{i+1}\rangle^2 - 16\langle\nabla_{U_{i+1}}U_{i+1},U_i\rangle^2
 + 16x_ix_{i+2} + 16x_{i+1}x_{i+2}\\[2ex]
    &=&-9b_i^2b_{i+2}^2d_{i+2}^{-2} - 9b_{i+1}^2b_{i+2}^2d_{i+2}^{-2}
    + 16x_ix_{i+2} + 16x_{i+1}x_{i+2}.
    \end{array}
\]
From the previous equations we get by a straightforward
calculation
\begin{equation}\label{thEqGauss}
\begin{array}{rcl}
0 & = & 8\lambda_i\lambda_{i+1}-2-12b_{i+2}^2
+9d_{i+2}^{-2}b_{i+2}^2(1-b_{i+2}^2) + 8(x_ix_{i+1} - x_ix_{i+2} - x_{i+1}x_{i+2})\\[2ex]
&& + 6d_{i+2}^{-1}\Big((b_{i+2}^2-b_{i+1}^2)x_i
-(b_{i+2}^2-b_i^2)x_{i+1} +b_{i+1}^2\lambda_i-b_i^2\lambda_{i+1}\Big).
\end{array}
\end{equation}

\medskip
The following proposition is the first crucial step towards the final classification.

\begin{proposition}\label{biconst}
The functions $b_1,b_2,b_3$ are constant.
\end{proposition}

\begin{proof}
Inserting $b_i^2,b_{i+1}^2,b_{i+2}^2$ according to (\ref{thbi}) into (\ref{thEqGauss})
gives
\begin{equation}\label{system}
\begin{array}{rcl}
0 & = & 8d_id_{i+1}d_{i+2}^{-2}(x_i+x_{i+1})^2 + 2d_i(4+2\lambda_id_{i+2}^{-1} - d_{i+2}^{-2})x_i\\[2ex]
&&- 2d_{i+1}(4-2\lambda_{i+1}d_{i+2}^{-1} - d_{i+2}^{-2})x_{i+1}
- 4(\lambda_i+\lambda_{i+1})x_{i+2} + d_{i+2}^{-2} - 2 + 4\lambda_i\lambda_{i+1}.
\end{array}
\end{equation}

This leads to a system of three quadratic equations of the form
\begin{eqnarray*}
(x_1+x_2)^2+\Lambda_{11}x_1+\Lambda_{12}x_2+\Lambda_{13}x_3&=&\Omega_1\\
(x_1+x_3)^2+\Lambda_{21}x_1+\Lambda_{22}x_2+\Lambda_{23}x_3&=&\Omega_2\\
(x_2+x_3)^2+\Lambda_{31}x_1+\Lambda_{32}x_2+\Lambda_{33}x_3&=&\Omega_3
\end{eqnarray*}
with some real constants $\Lambda_{ij}$ and $\Omega_i$.
We introduce new variables $y_1,y_2,y_3$ by means of $x_i=-y_i+y_{i+1}+y_{i+2}$.
Then the above system transforms into a system of three quadratic equations of the form
\[
4y_i^2+\bar{\Lambda}_{i1}y_1+\bar{\Lambda}_{i2}y_2+\bar{\Lambda}_{i3}y_3
=\Omega_i,\ i\in\{1,2,3\}.
\]
This system has only finitely many solutions (see e.g.\ Corollary 7 on page 233 in \cite{CLO92}).
It follows that the system (\ref{system}) has only finitely many solutions, and
as the coefficients of the system are constant, each solution must be constant.
Thus the functions $x_1,x_2,x_3$ are constant. From (\ref{thbi}) we see that
$b_1,b_2,b_3$ are constant. Finally,
\end{proof}

\begin{corollary}\label{thNo3Proj}
There exists a principal curvature $\lambda_i$ of $M$ such that the orthogonal
projection of $J\xi_p$ onto $T_{\lambda_i}(p)$ is equal to zero for all $p \in M$.
\end{corollary}

\begin{proof}
As $b_i$ is constant, equation (\ref{thdbii}) shows that $b_1b_2b_3 = 0$,
which implies the assertion.
\end{proof}

In view of the previous corollary we may assume $b_3=0$.
Moreover, if $b_1=0$ or $b_2=0$, then $M$ is a Hopf hypersurface.
In this case it follows from the classification of Hopf hypersurfaces with
constant principal curvatures in \cite{Be89} that $M$ is an open part of a tube of
radius $r\neq\ln(2+\sqrt{3})$ around a totally geodesic $\R
H^2\subset\CH$. Therefore we assume $b_1\neq 0$ and $b_2 \neq 0$ from now on.

Differentiating the constant function $b_2$ with respect to $U_3$ gives
\begin{equation}\label{x3}
x_3 = \lambda_3
\end{equation}
according to (\ref{thdbii1}), and differentiating $b_3 = 0$ with respect to $U_1$
and $U_2$ gives
\begin{equation}\label{x1x2}
4x_1 = 4\lambda_1 + 3d_2^{-1}b_1^2\ ,\ 4x_2 = 4\lambda_2 - 3d_1^{-1}b_2^2
\end{equation}
according to \eqref{thdbii1} and \eqref{thdbii2}. Inserting
these expressions for $x_1,x_2,x_3$ into the two equations for
$b_1^2$ and $b_2^2$ in (\ref{thbi}), and subtracting the two
resulting equations gives
\begin{equation}\label{b1b2lambda3}
d_1^{-1}b_2^2 - d_2^{-1}b_1^2 = 4\lambda_3.
\end{equation}
Together with the equation $b_1^2 + b_2^2 = 1$ this leads to
\begin{equation}\label{b1b2}
b_1^2 = d_2d_3^{-1}(4d_1\lambda_3 - 1)\ ,\ b_2^2 = -d_1d_3^{-1}(4d_2\lambda_3 + 1).
\end{equation}
We now insert the expressions for $x_1$ and $x_2$ according to (\ref{x1x2})
and the ones for $b_1^2$ and $b_2^2$ according to (\ref{b1b2}) into
the equation for $b_3^2 = 0$ in (\ref{thbi}), which gives
\[
2d_1\lambda_1 - 2d_2\lambda_2 + 6(d_2-d_1)\lambda_3 + 1 = 0.
\]
This equation is equivalent to
\[
(\lambda_1-\lambda_2)^2 - (\lambda_1+\lambda_2-4\lambda_3)^2 =
1-4\lambda_3^2.
\]
We now multiply equation (\ref{thEqGauss}) for $i=2$ with $d_1$
and equation (\ref{thEqGauss}) for $i=3$
with $d_2$, and then subtract the two resulting equations, which gives
\[
(10\lambda_1^2+10\lambda_2^2+6\lambda_1\lambda_2+1)\lambda_3
=2(\lambda_1+\lambda_2)(4\lambda_3^2+\lambda_1\lambda_2+1)
+6\lambda_3^3.
\]
We define $x=\lambda_1-\lambda_2$ and
$y=\lambda_1+\lambda_2-4\lambda_3$, which transforms the last two equations into
\[
x^2-y^2-1+4\lambda_3^2=0\ ,\ x^2(y+11\lambda_3) - y^3
+ \lambda_3 y^2+4(10\lambda_3^2-1)y
+2{\lambda_3}(34\lambda_3^2-7)=0.
\]
If $\lambda_3=0$ we immediately get
$\lambda_1,\lambda_2\in\{\pm1/2\}$. Hence, we assume
$\lambda_3\neq 0$. Inserting $x^2 = y^2+1-4\lambda_3^2$ into the
second equation gives a quadratic equation in $y$, which easily
leads to the possible solutions
\[
(x,y) = \left(\pm
\sqrt{1-3\lambda_3^2}\ ,\ -\lambda_3\right)\ \ {\rm and}\ \ (x,y)
= \left(\pm\frac{1}{4\lambda_3}\ ,\
\frac{1-8\lambda_3^2}{4\lambda_3}\right)\ ,
\]
where the first possibility only arises if $3\lambda_3^2 \leq 1$.
Taking into account that $\lambda_1$ and $\lambda_2$ are different
from $\lambda_3$, this eventually implies
\begin{equation}\label{eqLambda12}
\lambda_1 = \frac{1}{2}\left(3\lambda_3 -
\sqrt{1-3\lambda_3^2}\right)\ ,\ \lambda_2 =
\frac{1}{2}\left(3\lambda_3 + \sqrt{1-3\lambda_3^2}\right),
\end{equation}
where we assume without loss of generality that $\lambda_1 <
\lambda_2$. Obviously, we get a solution only if $3\lambda_3^2
\leq 1$. If $|\lambda_3|=1/2$ or $|\lambda_3|=1/\sqrt{3}$, then
the three principal curvatures cannot be different. Suppose  that
$1/2<|\lambda_3|<1/\sqrt{3}$. From (\ref{b1b2lambda3}) and
(\ref{eqLambda12}) we get
\[
\frac{b_1^2}{2\lambda_3(\lambda_3-\sqrt{1-3\lambda_3^2})}
+\frac{b_2^2}{2\lambda_3(\lambda_3+\sqrt{1-3\lambda_3^2})}=1.
\]
If $1/2<|\lambda_3|<1/\sqrt{3}$, elementary calculations show that
$0<2\lambda_3(\lambda_3-\sqrt{1-3\lambda_3^2})<1$ and
$0<2\lambda_3(\lambda_3+\sqrt{1-3\lambda_3^2})<1$. Therefore the
last equation is the equation of an ellipse centered at the origin
and with axes of length less than 1. Obviously such an ellipse has
no points of intersection with the circle $b_1^2+b_2^2=1$. This
shows that $|\lambda_3|<1/2$.

Therefore, we have proved

\begin{proposition}\label{thStep1}
Let $M$ be a connected real hypersurface in $\CH$ with three
distinct constant principal curvatures
$\lambda_1,\lambda_2,\lambda_3$, and assume that $M$ is not a
Hopf hypersurface. Then, with a suitable labelling of the principal
curvatures, we have
\[
2|\lambda_3| < 1\ ,\
2\lambda_1=3\lambda_3-\sqrt{1-3\lambda_3^2}\ ,\
2\lambda_2=3\lambda_3+\sqrt{1-3\lambda_3^2}.
\]
\end{proposition}

This result is the second crucial step towards the final classification.
We will now use Jacobi field theory to show that, under the assumptions of the
Proposition \ref{thStep1}, one of the equidistant hypersurfaces to $M$ is
an open part of the ruled minimal real hypersurface $W^3$.

For $r\in\R$ we define the smooth map $\Phi^r : M\to \CH,\ p
\mapsto \Phi^r(p)=\exp_p(r\xi_p)$, where $\exp_p$ is the
exponential map of $\CH$ at $p$. Geometrically this means that we
assign to $p$ the point in $\CH$ which is obtained by travelling
for the distance $r$ along the geodesic $c_p(t)=\exp_p(t\xi_p)$ in
direction of the normal vector $\xi_p$ (for $r >0$; for $r<0$ one
sets off in direction $-\xi_p$; and for $r=0$ there is no movement
at all). For $v \in T_pM$ we denote by $B_v$ the parallel vector
field along the geodesic $c_p$ with $B_v(0) = v$, and by $\zeta_v$
the Jacobi field along $c_p$ with $\zeta_v(0)=v$ and
$\zeta_v'(0)=-S_pv$. Note that $\zeta_v$ is the unique solution of
the linear differential equation
\[
4\zeta_v''-\zeta_v-3\langle \zeta_v,J\dot{c}_p\rangle
J\dot{c}_p=0\ ,\ \zeta_v(0)=v\ ,\ \zeta_v'(0)=-S_pv,
\]
where $\dot{c}_p$ denotes the tangent vector field of $c_p$ and
the prime $'$ indicates the covariant derivative of a vector field
along $c_p$. For $v \in T_{\lambda_i}(p)$ we have the explicit
expression
\[
\zeta_v(t)=f_i(t)B_v(t)+\langle v,J\xi\rangle g_i(t)J\dot{c}_p(t)
\]
with
\[
\begin{array}{rcl}
f_i(t)    &=& \displaystyle
        \cosh(t/2)-2\lambda_i\sinh(t/2),\\
\noalign{\medskip} g_i(t)    &=&
\displaystyle\left(\cosh(t/2)-1\right)
    \left(1+2\cosh(t/2)-2\lambda_i\sinh(t/2)\right).
\end{array}
\]
Finally, we define a vector field $\eta^r$ along the map $\Phi^r$
by $\eta^r_p=\dot{c}_p(r)$. The relation between the map $\Phi^r$,
the vector field $\eta^r$ and the Jacobi field $\zeta_v$ is given
by
\[
\zeta_v(r) = \Phi^r_*v\ ,\
\zeta_v'(r) = \enabla_{v}\eta^r,
\]
where $\Phi^r_*$ denotes the differential of $\Phi^r$.
We will see that there exists a real number $r \in \R$ such that
the map $\Phi^r$ has constant rank $3$ and the image is
locally a minimal real hypersurface ${\mathcal W}$. We then use the equation
$\zeta_v'(r) = \enabla_{v}\eta^r$ to obtain some information about
the second fundamental form of ${\mathcal W}$.

The following result was proved in \cite{BD06} for $n\geq
3$, but the same argument holds for $n=2$.

\begin{theorem}\label{ruledk}
Let $M$ be a $3$-dimensional connected submanifold in $\CH$.
Assume that there exists a unit vector field $Z$ tangent to the
maximal holomorphic subbundle of $TM \subset T\CH$ such that the
second fundamental form $I\!I$ of $M$ is given by the trivial
bilinear extension of $2I\!I(Z,J\xi) = \xi$ for all $\xi \in
\Gamma(\nu M)$. Then $M$ is holomorphically congruent to an open
part of the ruled minimal submanifold $W^3$.
\end{theorem}

If $\lambda_3= 0$, then $\lambda_1 = -1/2$ and $\lambda_2 = 1/2$, and
\eqref{b1b2} implies $2b_1^2 = 2b_2^2 = 1$. Thus, if we define
$Z = -b_1U_1 + b_2U_2$, we see from Theorem \ref{ruledk} that $M$ is
holomorphically congruent to an open part of the ruled minimal
hypersurface $W^3$.

If $\lambda_3 \neq 0$, then $0 < 2|\lambda_3| < 1$
by Proposition \ref{thStep1}. Thus we can write $2\lambda_3 = \tanh(r/2)$
with some nonzero real number $r$. Using the equation
$\Phi^r_*v = \zeta_v(r)$ and the explicit expression for the
Jacobi fields, we obtain
\[
\Phi^r_* U_3(p) = \sech(r/2)B_{U_3(p)}(r)
\]
and
\[
\left( \begin{array}{@{}c@{}} \Phi^r_*U_1(p) \\ \Phi^r_*U_2(p)
\end{array} \right) = D(r) \left( \begin{array}{@{}c@{}} B_{U_1(p)}(r) \\
B_{U_2(p)}(r)
\end{array} \right)
\]
with
\[
D(t) = \left(
\begin{array}{@{}cc@{}}
f_1(t)+b_1^2 g_1(t)   &   b_1b_2g_1(t)\\
b_1b_2g_2(t)  & f_2(t)+b_2^2 g_2(t)
\end{array}\right)
\]
for all $p \in M$. As $\det(D(r)) = \sech^3(r/2)$, we can now conclude that
$\Phi^r_*$ has maximal rank everywhere. This means that for every
point in $M$ there exists an open neighborhood $\mathcal{V}$
around that point such
that $\mathcal{W}=\Phi^r(\mathcal{V})$ is an embedded real
hypersurface of $\CH$ and $\Phi^r:\mathcal{V}\to\mathcal{W}$ is a
diffeomorphism. Let ${\mathcal V}$ be such an open neighborhood. We fix a point
$p \in {\mathcal V}$ and define $q = \Phi^r(p) \in
{\mathcal W}$. The tangent space $T_q{\mathcal W}$ of ${\mathcal
W}$ at $q$ is obtained by parallel translation of $T_p{\mathcal
V}$ along the geodesic $c_p$ from $p = c_p(0)$ to $q = c_p(r)$,
and $\eta^r_p$ is a unit normal vector of ${\mathcal W}$ at $q$.

For the shape operator $S^r$ of ${\mathcal W}$ we have
$S^r_{\eta^r_p}\Phi^r_*v = -\enabla_{v}\eta^r = -\zeta_v'(r)$
for all $v \in T_pM$.
Since $f_3'(r) = 0$, we immediately get
\begin{equation}\label{eqSBU3}
S^r_{\eta^r_p}B_{U_3(p)}(r) = 0,
\end{equation}
and by putting $v = U_1(p)$ and $v=U_2(p)$ we get
\[
\left( \begin{array}{@{}c@{}} S^r_{\eta^r_p}B_{U_1(p)}(r) \\
S^r_{\eta^r_p}B_{U_2(p)}(r) \end{array}\right)
= C(r) \left( \begin{array}{@{}c@{}} B_{U_1(p)}(r) \\
B_{U_2(p)}(r)
\end{array} \right)
\]
with $C(r) = -D'(r)D(r)^{-1}$. A tedious calculation shows that
$\det(D'(r)) = -\sech^3(r/2)/4$ and $(\det(D))'(r) = 0$, which
implies
\[
\det(C(r)) = \frac{\det(D'(r))}{\det(D(r))} = -\frac{1}{4} \quad
\mbox{and}\quad  \tr(C(r)) = -\frac{(\det(D))'(r)}{\det(D(r))} =
0.
\]
From this we easily see that the eigenvalues of $C(r)$ are $\pm
1/2$. Altogether we now get that ${\mathcal W}$ has three distinct
constant principal curvatures $0$, $+1/2$ and $-1/2$.
According to the classification of real Hopf hypersurfaces in
\cite{Be89} we see that $\mathcal{W}$ cannot be a Hopf
hypersurface. As $U_3(p)$ belongs to the maximal holomorphic subspace
of $T_pM$ and parallel translation in $\C H^2$
along the geodesic $c_p$ commutes with the
complex structure $J$, the vector $B_{U_3(p)}(r)$ belongs to the maximal
holomorphic subspace of $T_q{\mathcal W}$. As we have seen above,
$B_{U_3(p)}(r)$ is a principal curvature vector corresponding to
the principal curvature $0$.
It follows that we can write $J\eta^r_p$ as a linear combination of
principal curvature vectors corresponding to the principal curvatures
$\pm 1/2$. Applying \eqref{b1b2} to ${\mathcal W}$ instead of $M$ we see
that we can write $J\eta^r = (X_+ + X_-)/\sqrt{2}$ with suitable unit
vector fields $X_+$ and $X_-$ corresponding to
the principal curvatures $+1/2$ and $-1/2$, respectively.
Defining $Z=(X_+-X_-)/\sqrt{2}$ we get that the second fundamental form of
$\mathcal{W}$ satisfies the formula of Theorem \ref{ruledk}.
It follows that ${\mathcal W}$ is holomorphically congruent to an open part
of the ruled real hypersurface $W^3$. From this we eventually
conclude that $M$ is holomorphically congruent to an open part of
an equidistant hypersurface to $W^3$.

This finishes the proof of Theorem \ref{mainresult}.


\noindent {\sc Department of Mathematics, University College,
Cork, Ireland}

\smallskip
\noindent {\sc Department of Geometry and Topology, Faculty of
Mathematics, University of Santiago de Compostela, Spain}

\end{document}